\documentclass{article}
\usepackage{calc}
\usepackage{tikz-cd}
\usepackage{tensor}
\usepackage{csquotes}
\usepackage{fancyvrb}
\usepackage{verbatim}
\usepackage{amsmath}
\usepackage{microtype}
\usepackage{amssymb}
\usepackage{amsthm,thmtools}
\usepackage{nameref,hyperref}
\usepackage{cleveref}
\declaretheorem[numberwithin=section]{theorem}

\declaretheorem[sibling=theorem]{non-example}
\declaretheorem[sibling=theorem]{lemma}

\declaretheorem[sibling=theorem]{proposition}
\declaretheorem[style=definition,sibling=theorem]{definition}

\declaretheorem[style=definition,sibling=theorem]{notation}

\declaretheorem[style=remark,sibling=theorem]{remark}
\title{The Algebroid of a Groupoid in a Tangent Category}
\author{Matthew Burke\footnote{The author acknowledges the support of a Postdoctoral Scholarship from the Department of Mathematics and Statistics at the University of Calgary. This preprint is part of a larger project that is joint work with Robin Cockett and Geoff Cruttwell.}}
\begin{document}
\maketitle
\begin{abstract}
  We generalise the construction of the Lie algebroid of a Lie groupoid so that it can be carried out in any tangent category. First we reconstruct the bijection between left invariant vector fields and source constant tangent vectors based at an identity element for a groupoid in a category equipped with an endofunctor that has a retraction onto the identity functor. Second we use the full structure of a tangent category to construct the algebroid of a groupoid. Finally we show how the classical result concerning the splitting of the tangent bundle of a Lie group can be carried out for any pregroupoid.
\end{abstract}
\tableofcontents

\section{Introduction}
\label{sec:introduction}

Recall that in classical Lie theory we approximate a Lie group using an algebraic structure called a Lie algebra. To construct the Lie algebra that approximates a particular Lie group $\mathbb G = (G,\mu,e)$ we first consider the tangent space of the underlying smooth manifold $G$ at the identity element $e$ of the underlying group. To define the Lie bracket $[v,w]$ of two tangent vectors $v$ and $w$ based at the identity element we use the group multiplication $\mu$ to extend $v$ and $w$ to (left invariant) vector fields $X$ and $Y$ on $G$ then define $[v,w]$ to be the evaluation of $[X,Y]$ at the identity element. The Lie bracket defined in this way is well defined due to two facts. Firstly there is a bijection between tangent vectors at the identity and left invariant vector fields. Secondly if $X$ and $Y$ are left invariant vector fields then $[X,Y]$ is left invariant also. For more details concerning the classical theory of Lie groups and Lie algebras please see \cite{MR1738431}. A more conceptual way of phrasing the correspondence between tangent vectors at the identity element and left invariant vector fields is to say that the tangent space $TG$ splits into the product $G\times T_e G$ where $TG$ denotes the tangent space of $G$ and $T_e G$ the subspace of $TG$ consisting of vectors based at the identity element.

In the established generalisation of Lie theory involving Lie groupoids and Lie algebroids (see for instance \cite{MR1973056}) we can also construct a Lie algebroid that approximates a Lie groupoid. Recall that a \emph{Lie groupoid} is a groupoid in the category $Man$ of smooth manifolds such that the source and target maps are submersions. Recall further that a \emph{Lie algebroid} is a smooth vector bundle $A \rightarrow M$ together with a bundle homomorphism $\rho\colon A \rightarrow TM$ such that the space of sections $\Gamma(A)$ is a Lie algebra satisfying the following Leibniz law: for all $X,Y\in \Gamma(A)$ and $f\in C^{\infty}(M)$ the equality
\begin{gather*}[X,fY]=\rho(X)(f)\cdot Y+f[X,Y]\end{gather*}
holds. The construction of the Lie algebroid that is the linear approximation to a Lie groupoid is analogous to the construction of the Lie algebra that approximates a Lie group described above. First we consider the space $T^s_M G$ of tangent vectors that are based at an identity element of $\mathbb G$ and are parallel to the source fibres of $\mathbb G$. The reason we insist on vectors that are parallel to the source fibres of $\mathbb G$ is that we can then use the composition of the Lie groupoid to extend a section of the natural projection $\pi^s_M: T^s_M G \rightarrow M$ to a vector field on $G$. Then we can use the standard Lie bracket of vector fields to induce a Lie bracket on sections of $\pi^s_M$. For more details on the approximation of a Lie groupoid by a Lie algebroid please see Section 3.5 of \cite{MR2157566}.

In this paper we show how to generalise this technique of approximating a groupoid by an algebroid to any category $\mathbb X$ equipped with a tangent structure (or {\em tangent category}) as defined in Section 2 of \cite{MR3192082}. In fact some of the theory can be carried out at a higher level of generality. In \Cref{sec:invariant-vector-fields} we show that the bijection between left invariant vector fields and sections of $\pi^s_M$ holds in any category equipped with an endofunctor that has a retraction onto the identity functor. One example of such an endofunctor that is not necessarily a tangent structure is the diagonal functor $T=\Delta: Set_* \rightarrow Set_*$ on the category of pointed sets defined by $X \mapsto X\times X$. Then we can take $p$ to be the second projection and $0:X \rightarrow X\times X$ to be defined by $x\mapsto (*,x)$. Another example uses the double negation operation in a topos. Suppose that $\mathcal E$ is a topos and that $\neg\neg : \mathcal E \rightarrow \mathcal E$ is the double negation endofunctor. Then we can define $T=\Delta(\neg\neg)$ which takes $X\mapsto \{(x,y)\in X\times X: \neg\neg(x=y)\}$. Then we can define $p$ to be the second projection and $0$ to be the diagonal. In the well-adapted model of synthetic differential geometry called the Dubuc topos \cite{MR1055025} this $T$ defines the bundle of germs on a manifold.

In \Cref{sec:the-algebroid-of-a-groupoid-in-a-tangent-category} we assume the full structure of a tangent category to construct the algebroid that is the linear approximation to a groupoid in $\mathbb X$. In \Cref{sec:torsoids-in-tangent-categories} we show that the splitting of the tangent bundle holds not only for groupoids but also for pregroupoids (as defined in \cite{math/0502075}) equipped with an extra arrow satisfying a certain equation. In this section also we do not need the full structure provided by a tangent category and work in a category $\mathbb X$ equipped with a functor that has a retraction onto the identity functor.

\section{Invariant Vector Fields}\label{sec:invariant-vector-fields}

In this section we construct an abstract bijection which generalises the classical bijection between left invariant vector fields on the arrow space of a Lie groupoid and sections of the bundle of source constant tangent vectors based at an identity element of the groupoid. These constructions only require part of the structure of a tangent category. Therefore in this section we assume only the existence of: a category $\mathbb X$, an endofunctor $T:\mathbb X \rightarrow \mathbb X$, a natural transformation $p:T \Rightarrow 1_{\mathbb X}$ and a natural transformation $0:1_{\mathbb X} \Rightarrow T$. We further assume that $0$ is a section of $p$. In addition the definitions we make in this section require that certain pullbacks exist in $\mathbb X$. Therefore we assume without further comment that these pullbacks exist and moreover that they are preserved by $T$.

\subsection{Basic Definitions}

We define the arrows $\pi^s$ and $\pi^s_M$ of $\mathbb X$ that generalise the two classical bundles that we are interested in. The arrow $\pi^s$ corresponds to the bundle of source constant tangent vectors. The arrow $\pi^s_M$ corresponds to the bundle of source constant tangent vectors based at an identity element.

\begin{notation}
  Let $\mathbb{X}$ be a category and $T:\mathbb X \rightarrow \mathbb X$ a functor. Let $p: T \Rightarrow  1_{\mathbb X}$ and $0:1_{\mathbb X} \Rightarrow  T$ be natural transformations such that $0$ is a section of $p$. Let $\mathbb{G}=G \rightrightarrows M$ be a groupoid in $\mathbb{X}$ with source, target, inverse and composition denoted by $s$, $t$, $(-)^{-1}$ and $\mu$ respectively. This means that the pullback $G \tensor[_t]{\times}{_s} G$ exists and further we assume that $T$ preserves this pullback.
\end{notation}

\begin{definition}\label{pullbacks-of-the-tangent-bundle}
  We denote by $\pi^{s}:T^{s} G \rightarrow M$ the pullback
  \begin{equation}
    \begin{tikzcd}
      T^{s} G \rar{\iota^{s}} \dar{\pi^{s}} & TG \dar{Ts}\\
      M \rar{0} & TM
    \end{tikzcd}
  \end{equation}
  consisting of $s$-constant tangent vectors of $G$ and write $p^s = \iota^s p$. Finally we denote by $\pi^{s}_M:T^{s}_M G \rightarrow M$ the pullback
  \begin{equation}
    \begin{tikzcd}
      T^{s}_M G \rar{\iota_M} \dar{\pi^{s}_M} & T^{s} G \dar{p^s}\\
      M \rar{e} & G
    \end{tikzcd}
  \end{equation}
  and write $\iota^s _M = \iota_M \iota^s$ and $p^s _M = \iota_M p^s$.
\end{definition}

\begin{remark}
  The arrow $\pi^{s}_M:T^{s}_M G \rightarrow M$ is equivalently defined by the double pullback
  \begin{equation}\label{double-pullback-definition}
    \begin{tikzcd}
      T^{s}_M G \arrow{dr}{\iota^{s}_M} \arrow{dd}{\pi^{s}_M} \arrow{rr}{q} & {} & M \dar{0}\\
      {} & TG \dar{p} \rar{Ts} & TM\\
      M \rar{e} & G & {}
    \end{tikzcd}
  \end{equation}
\end{remark}

\begin{remark}
  In fact $q=\pi^{s}_M$. Indeed if
  \begin{equation}\label{cone-over-double-pullback}
    \begin{tikzcd}
      Z \arrow{dr}{z_1} \arrow{dd}{z_0} \arrow{rr}{z_2} & {} & M \dar{0}\\
      {} & TG \dar{p} \rar{Ts} & TM\\
      M \rar{e} & G & {}
    \end{tikzcd}
  \end{equation}
  is an arbitrary cone over the diagram defining $p^{s}_M$ then we can paste the commutative square
  \begin{equation}
    \begin{tikzcd}
      TG \rar{Ts} \dar{p} & TM \dar{p}\\
      G \rar{s} & M
    \end{tikzcd}
  \end{equation}
  onto the bottom right of \eqref{cone-over-double-pullback} and conclude that $z_2=z_20p=z_0 e s=z_0$.

  In fact we can make a further simplification in the definition of $T^{s}_M G$. Since the components $z_0$ and $z_2$ of any cone over the limit defining $T^{s}_M G$ are equal we can conflate the two copies of $M$ in \eqref{cone-over-double-pullback} and conclude that $T^{s}_M G$ is the limit of the diagram
  \begin{equation}
    \begin{tikzcd}
      TG \rar{Ts} \dar{p} & TM \\
      G & M \lar{e} \uar{0}
    \end{tikzcd}
  \end{equation}
  and so in fact $\pi^{s}_M$ is equivalently defined using the pullback
  \begin{equation}
    \begin{tikzcd}
      T^{s}_M G \rar{s} \dar{\pi^{s}_M} & TG \dar{(Ts,p)}\\
      M \rar{(0,e)} & TM \times G
    \end{tikzcd}
  \end{equation}
\end{remark}

\subsection{Extending and Restricting Vector Fields}\label{sec:bijection}

In this section we recall how to postcompose an element of a groupoid with a source constant tangent vector and then how to define a left invariant vector field. With these definitions in place we show how to extend a section of $\pi^s _M$ to a left invariant vector field and restrict a left invariant vector field to a section of $\pi^s_M$. Finally we show that these extension and restriction functions are inverses.

\begin{notation}
  The arrow $(\pi_0 0,\pi_1 \iota^s)T\mu:G \tensor[_t]{\times}{_{\pi^s}} T^s G \rightarrow TG$ factors through $T^s G$ because 
  \begin{align}
    (\pi_0 0,\pi_1 \iota^s)T\mu T s & = (\pi_0 0,\pi_1 \iota^s) \pi_0 Ts\\
    & = \pi_0 0 Ts\\
    & = \pi_0 s 0
  \end{align}
  We call this unique factorisation $\mu^s:G \tensor[_t]{\times}{_{\pi^s}} T^s G \rightarrow T^s G$.
  Let $\Gamma(M,T^s _M G)$ be the set of sections of $\pi^s _M$.
\end{notation}

\begin{definition}
	A \emph{left invariant vector field $X:G \rightarrow T^s G$} is a section of $p^s:T^s G \rightarrow G$ such that
	\begin{equation}
		\begin{tikzcd}
			G\tensor[]{\times}{}T^s G \rar{\mu^s} & T^s G\\
			G\tensor[]{\times}{} G \rar{\mu} \uar{(\pi_0,\pi_1 X)} & G \uar{X}
		\end{tikzcd}
	\end{equation}
  commutes.
  We write ${\sf LeftInv}(G,T^s G)$ for the set of sections of $p^s$ that are left invariant.
\end{definition}

\begin{lemma}
  \begin{align}
    \mu^s p^s &= \mu^s \iota^s p\\
    & = (\pi_0 0,\pi_1 \iota^s)T\mu p\\
    & = (\pi_0 0,\pi_1 \iota^s) p \mu\\
    & = (\pi_0,\pi_1 p^s) \mu
  \end{align}
\end{lemma}

\begin{lemma}
  Every section $v:M \rightarrow T^s _M G$ of $\pi^s _M:T^s _M G \rightarrow M$ can be extended to a left invariant vector field $v^{\wedge}:G \rightarrow T^s G$ of $p^s:T^s G \rightarrow G$.
  \begin{proof}
    Define $v^{\wedge}$ as the composite
    \begin{equation}
      G \xrightarrow{(1_G,t)} G \tensor[_t]{\times}{_{1_M}} M \xrightarrow{(\pi_0,\pi_1 v \iota_M)} G \tensor[_t]{\times}{_{\pi^s}} T^s G \xrightarrow{\mu^s} T^s G
    \end{equation}
    Then $v^{\wedge}$ is a section of $p^s$ because
    \begin{align}
      v^{\wedge} p^s & = (1_G, tv\iota_M)\mu^s p^s \\
      & = (1_G, tv\iota_M)(\pi_0,\pi_1 p^s) \mu\\
      & = (1_G,tv\iota_M p^s) \mu\\
      & = (1_G, t v \pi^s _M e)\mu\\
      & = (1_G, te)\mu = 1_G
    \end{align}
    Finally $v^{\wedge}$ is left invariant because
    \begin{align}
      (\pi_0,\pi_1 v^{\wedge}) \mu^s \iota^s &= (\pi_0, \pi_1 (1_G, tv\iota_M)\mu^s) (\pi_0 0,\pi_1 \iota^s)T\mu\\
      & = (\pi_0 0, \pi_1 (1_G, tv\iota_M)\mu^s \iota^s) T\mu\\
      & = (\pi_0 0, \pi_1 (1_G, tv\iota_M)(\pi_0 0,\pi_1 \iota^s)T\mu) T\mu\\
      & = (\pi_0 0, \pi_1 (1_G 0,tv\iota_M \iota^s)T\mu) T\mu\\
      & = (\pi_0 0, (\pi_1 0,\pi_1 tv\iota_M \iota^s)T\mu) T\mu\\
      & = ((\pi_0 0,\pi_1 0)T\mu,\pi_1 tv\iota_M \iota^s) T\mu\\
      & = ((\pi_0,\pi_1)\mu 0,\pi_1 tv\iota_M \iota^s) T\mu\\
      & = ((\pi_0,\pi_1)\mu,\pi_1 tv\iota_M)(\pi_0 0,\pi_1 \iota^s) T\mu\\
      & = ((\pi_0,\pi_1)\mu,\pi_1 tv\iota_M)\mu^s \iota^s\\
      & = ((\pi_0,\pi_1)\mu,(\pi_0,\pi_1) tv\iota_M)\mu^s \iota^s\\
      & = (\pi_0,\pi_1)\mu(1_G,tv\iota_M)\mu^s \iota^s\\
      & = \mu v^{\wedge}\iota^s
    \end{align}
    and so $(\pi_0,\pi_1 v^{\wedge}) \mu^s = \mu v^{\wedge}$ because $\iota^s$ is a monomorphism.
  \end{proof}
\end{lemma}

\begin{lemma}
  Every section $X:G \rightarrow T^s G$ of $p^s:T^s G \rightarrow G$ restricts to a section $X^{\vee}:M \rightarrow T^s _M G$ of $\pi^s _M: T^s _M G \rightarrow M$.
  \begin{proof}
    The arrow $e X: M \rightarrow T^s G$ factors through $T^s _M G$ because $e X p^s = e$.
    Let $X^{\vee}$ be this unique factorisation induced by $(1_M,e X)$.
    This means that $X^{\vee}\pi^s _M = 1_M$ and so $X^{\vee}$ is a section of $\pi^s _M$.
  \end{proof}
\end{lemma}

\begin{proposition}\label{bijection-of-sections}
  The function $(-)^{\vee}:\Gamma(M,T^s _M G) \rightarrow {\sf LeftInv}(G,T^s G)$ has inverse $(-)^{\wedge}$.
  \begin{proof}
    Let $v:M \rightarrow T^s _M G$ be a section of $\pi^s _M: T^s_M G \rightarrow G$.
    Then
    \begin{align}
      (v^{\wedge})^{\vee} & = ((1_G, tv\iota_M)\mu^s)^{\vee}\\
      & = (1_M, e (1_G, tv\iota_M)\mu^s)\\
      & = (1_M, (e, v\iota_M)\mu^s)\\
      & = (1_M, v\iota_M)\\
      & = v
    \end{align}
    Let $X:G \rightarrow T^s G$ be a left invariant vector field.
    Then
    \begin{align}
      (X^{\vee})^{\wedge} & = (1_M,e X)^{\wedge}\\
      & = (1_G, t(1_M,e X)\iota_M)\mu^s\\
      & = (1_G, te X)\mu^s\\
      & = (1_G,te)(\pi_0,\pi_1 X)\mu^s\\
      & = (1_G, te) \mu X = X
    \end{align}
    where the penultimate equality is due to the fact that $X$ is left invariant.
  \end{proof}
\end{proposition}

\section{The Algebroid of a Groupoid in a Tangent Category}\label{sec:the-algebroid-of-a-groupoid-in-a-tangent-category}

In this section we generalise the classical construction that assigns a Lie algebroid to each Lie groupoid. To carry out this work we must use the full structure of a tangent category. Therefore for this section only we assume that the category $\mathbb X$ that we are working in is a tangent category. If $\mathbb G = G \rightrightarrows M$ is a groupoid in $\mathbb X$ then the underlying bundle of the algebroid will be the bundle $\pi^s_M$ defined in \Cref{pullbacks-of-the-tangent-bundle}. We define the Lie bracket of two sections of $\pi^s_M$ in two steps. First we use the bijection of \Cref{bijection-of-sections} to obtain two left invariant vector fields on $G$. Second we use the description of the Lie bracket of vector fields in \cite{MR3192082} to define a Lie bracket. Then in order to prove that this bracket is well defined we check that the vector field thus obtained is source constant and left invariant.

\begin{notation}
  In this section we make further assumptions about $\mathbb X$, $T$, $p$ and $0$ so that they form part of a tangent structure as defined in \cite{MR3192082}. This means that we have the natural transformations $c:T^2 \rightarrow  T^2$ and $\ell: T \rightarrow  T^2$, and that $p:T \Rightarrow 1_{\mathbb X}$  is an additive bundle.
\end{notation}

The following result is Lemma 2.13 in \cite{MR3192082}.

\begin{lemma}
  In any tangent category the following is a (triple) equaliser diagram:
  \begin{equation}\label{triple-equaliser}
    \begin{tikzcd}
    TG \rar{l} & T^2 G \rar[][description]{p} \rar[yshift=+2ex]{Tp} \rar[yshift=-2ex][swap]{pp0} & TG
  \end{tikzcd}
  \end{equation}
\end{lemma}

The following is Definition 3.14 in \cite{MR3192082}.

\begin{definition}
  In a tangent category with vector fields $X,Y:G \rightarrow TG$ the \emph{Lie bracket} of $X$ and $Y$ is defined to be the morphism
  \begin{equation}
    [X,Y] := \{XTY - YTX c\}: G \rightarrow TG
  \end{equation}
  where $c$ is the canonical flip and $\{f\}$ denotes the factorisation of $f$ through $l$ for any arrow $f$ that equalises diagram \eqref{triple-equaliser}.
\end{definition}

\begin{lemma}\label{lie-bracket-is-hom}
  If $X,Y\in {\sf LeftInv}(G,T^s G)$ then $[X\iota^s,Y\iota^s]$ factors through $T^s G$.
  (This means that there is an arrow $[X,Y]:G \rightarrow T^s G$ such that $[X,Y]\iota^s = [X\iota^s,Y\iota^s]$.)
  \begin{proof}
    First
    \begin{equation}
      \begin{tikzcd}
        G \rar{X} \dar{s} & T^s G \rar{\iota^s} \dar{\pi^s} & TG \rar{TY} \dar{Ts} & T(T^s G) \rar{T(\iota^s)} \dar{T(\pi^s)} & T^2 G \dar{T^2 s}\\
        M \rar{=} & M \rar{0} & TM \rar{=} & TM \rar{T0} & T^2 M
      \end{tikzcd}
    \end{equation}
    where the first and third squares commute because $\pi^s = \pi^s 0 p = \iota^s Ts p = \iota^s p s = p^s s$.
    Similarly $Y \iota^s T(X\iota^s) T^2 s = s 0 T0$.
    Now we show that the arrow $[X\iota^s,Y\iota^s]:G \rightarrow TG$ factors through $T^s G$:
    \begin{align}
      [X\iota^s,Y\iota^s] Ts l & = \{X\iota^sT(Y\iota^s) - Y\iota^sT(X\iota^s) c\} l T^2 s\\
      & = (X\iota^sT(Y\iota^s) - Y\iota^sT(X\iota^s) c) T^2 s\\
      & = X\iota^sT(Y\iota^s)T^2 s - Y\iota^sT(X\iota^s)c T^2 s\\
      & = X\iota^sT(Y\iota^s)T^2s - Y\iota^sT(X\iota^s)T^2 s c\\
      & = s0 T 0 - s 0 T 0 c\\
      & = s 0 T0\\
      & = s 0 l
    \end{align}
    and so $[X\iota^s,Y\iota^s] Ts = (s-s) 0$ because $l$ is a monomorphism.
  \end{proof}
\end{lemma}

\begin{lemma}
  If $X,Y\in {\sf LeftInv}(G,T^s G)$ then the arrow $[X,Y]$ defined in \Cref{lie-bracket-is-hom} is in $\Gamma(G,T^s G)$.
  \begin{proof}
    We show that $[X\iota^s,Y\iota^s]$ is a section of $p$:
    \begin{align}
        [X\iota^s,Y\iota^s] p & = \{X\iota^sT(Y\iota^s) - Y\iota^sT(X\iota^s) c\} p\\
        & = (X\iota^sT(Y\iota^s) - Y\iota^sT(X\iota^s) c) pp\\
        & = X\iota^sT(Y\iota^s) pp\\
        & = (X\iota^s p)(Y\iota^s p) \\
        & = 1_G 1_G = 1_G
      \end{align}
      as required.
  \end{proof}
\end{lemma}

\begin{lemma}
  If $X,Y\in {\sf LeftInv}(G,T^s G)$ then the arrow $[X,Y]$ defined in \Cref{lie-bracket-is-hom} is in ${\sf LeftInv}(G,T^s G)$.
  \begin{proof}
      \begin{align*}
        (\pi_0,\pi_1 [X,Y]) \mu^s \iota^s l &= (\pi_0,\pi_1 [X,Y]) (\pi_0 0,\pi_1 \iota^s)T\mu l\\
        & = (\pi_0 0, \pi_1 [X,Y] \iota^s)T\mu l\\
        & = (\pi_0 0, \pi_1 [X\iota^s,Y\iota^s])T\mu l\\
        & = (\pi_0 0, \pi_1 [X\iota^s,Y\iota^s]) l T^2 \mu\\
        & = (\pi_0 0 l, \pi_1 [X\iota^s,Y\iota^s]l) T^2 \mu\\
        & = (\pi_0 0 l, \pi_1 (X\iota^sT(Y\iota^s) - Y\iota^sT(X\iota^s) c)) T^2 \mu\\
        & = (\pi_0 0 l, \pi_1 X\iota^sT(Y\iota^s))T^2 \mu - (\pi_0 0 l,\pi_1Y\iota^sT(X\iota^s) c)T^2 \mu \\
        & = (\pi_0 0 T0,\pi_1 X \iota^s T(Y \iota^s))T^2\mu - (\pi_0 0 T0,\pi_1 Y \iota^s T(X\iota^s))T^2\mu c\\
        & = (\pi_0 0,\pi_1 X \iota^s)(T\pi_0 T0,T\pi_1 T(Y \iota^s))T^2\mu - (\pi_0 0,\pi_1 Y \iota^s)(T\pi_0 T0,T\pi_1 T(X \iota^s))T^2\mu c\\
        & = (\pi_0 0,\pi_1 X \iota^s)T(\pi_0 0,\pi_1 Y \iota^s)T^2\mu - (\pi_0 0,\pi_1 Y \iota^s)T(\pi_0 0,\pi_1 X \iota^s)T^2\mu c\\
        & = (\pi_0 0,\pi_1 X \iota^s)T((\pi_0,\pi_1 Y)(\pi_0 0,\pi_1 \iota^s)T\mu) - (\pi_0 0,\pi_1 Y \iota^s)T((\pi_0,\pi_1 X)(\pi_0 0,\pi_1 \iota^s)T\mu) c\\
        & = (\pi_0 0,\pi_1 X \iota^s)T((\pi_0,\pi_1 Y)\mu^s\iota^s) - (\pi_0 0,\pi_1 Y \iota^s)T((\pi_0,\pi_1 X)\mu^s\iota^s) c\\
        & = (\pi_0 0,\pi_1 X \iota^s)T(\mu Y\iota^s) - (\pi_0 0,\pi_1 Y \iota^s)T(\mu X\iota^s) c\\
        & = (\pi_0,\pi_1 X)(\pi_0 0,\pi_1 \iota^s)T\mu T(Y\iota^s) - (\pi_0,\pi_1 Y)(\pi_0 0,\pi_1 \iota^s)T\mu T(X\iota^s) c\\
        & = (\pi_0,\pi_1 X)\mu^s\iota^sT(Y\iota^s) - (\pi_0,\pi_1 Y)\mu^s\iota^sT(X\iota^s) c\\
        & = \mu X\iota^sT(Y\iota^s) - \mu Y\iota^sT(X\iota^s) c\\
        & = \mu (X\iota^sT(Y\iota^s) - Y\iota^sT(X\iota^s) c)\\
        & = \mu [X\iota^s,Y\iota^s] l\\
        & = \mu [X,Y] \iota^s l
      \end{align*}
      and so $(\pi_0,\pi_1 [X,Y]) \mu^s = \mu [X,Y]$ because $\iota^s l$ is a monomorphism.
  \end{proof}
\end{lemma}

\section{Splitting the Tangent Bundle of a Pregroupoid}
\label{sec:torsoids-in-tangent-categories}

In this section we generalise the classical result that the bundle of source constant tangent vectors of a Lie groupoid splits as the product of the target bundle of the Lie groupoid and the bundle of source constant tangent vectors to the Lie groupoid based at an identity element. In fact we generalise this result in two directions. Firstly we replace the category of smooth manifolds with an arbitrary category equipped with a endofunctor which has a retraction onto the identity functor. Secondly we prove our results for the generalisation of groupoids called pregroupoids defined in \cite{math/0502075}.

\subsection{Pregroupoids and an Involution}

In this section we recall the definition of a pregroupoid. In addition we construct from this structure an involution that will be important in the sequel. Explicitly in \Cref{sec:splitting-of-the-vertical-bundle} we assume the existence of an extra arrow in addition to the pregroupoid structure; in the presence of this extra arrow the involution we define in this section will pullback in a natural way to give the splitting of the tangent bundle.

\begin{notation}
  Let $\mathbb{X}$ be a category and $T: \mathbb X \rightarrow \mathbb X$ a functor. Let $p:T \Rightarrow 1_{\mathbb X}$ and $0: 1_{\mathbb X} \Rightarrow T$ be natural transformations such that $0$ is a section of $p$.
\end{notation}

The following definition is in Section 1 of \cite{math/0502075}.

\begin{definition}
  A \emph{pregroupoid in $\mathbb{X}$} consists of a span
  \begin{equation}
    \begin{tikzcd}
      & X \drar{\beta} \dlar[swap]{\alpha} & \\
      A & & B
    \end{tikzcd}
  \end{equation}
  in $\mathbb{X}$ and a pre-composition $\mu=-\circ_{-}-:X\tensor[_{\beta}]{\times}{_{\beta}}X\tensor[_{\alpha}]{\times}{_{\alpha}}X \rightarrow X$ satisfying the following equations:-
  \begin{itemize}
    \item $\alpha(x\circ_y z) = \alpha(x)$;
    \item $\beta(x\circ_y z) = \beta(z)$;
    \item if $\beta(x) = \beta(y)$ then $x\circ_y y = x$;
    \item if $\alpha(x) = \alpha(y)$ then $x\circ_x y = y$;
    \item if $\beta(x) = \beta(y) = \beta(z)$ and $\alpha(z) = \alpha(w)$ then $x\circ_y(y\circ_z w) = x\circ_z w$;
    \item if $\beta(x) = \beta(y)$ and $\alpha(y) = \alpha(z) = \alpha(w)$ then $(x\circ_y z)\circ_z w = x\circ_y w$.
  \end{itemize}
\end{definition}

\begin{lemma}
  The arrow
  $$X\tensor[_{\beta}]{\times}{_{\beta}}X\tensor[_{\alpha}]{\times}{_{\alpha}}X \xrightarrow{\check{\xi}} X\tensor[_{\beta}]{\times}{_{\beta}}X\tensor[_{\alpha}]{\times}{_{\alpha}}X$$
  defined by $(\pi_1,\pi_0,\pi_0\circ_{\pi_1}\pi_2)$ is an involution.
  \begin{proof}
    \begin{align*}
      \check{\xi}\check{\xi} & = (\pi_1, \pi_0, \pi_0 \circ_{\pi_1}\pi_2)(\pi_1, \pi_0, \pi_0 \circ_{\pi_1}\pi_2)\\
                         & = (\pi_0,\pi_1,\pi_1\circ_{\pi_0}(\pi_0\circ_{\pi_1}\pi_2))\\
      & = (\pi_0,\pi_1,\pi_2) = 1
    \end{align*}
  \end{proof}
\end{lemma}

\begin{definition}
  The arrow \emph{$\pi^{\alpha}$} is defined as the pullback
  \begin{equation}
    \begin{tikzcd}
      T^{\alpha}X \dar{\pi^{\alpha}} \rar{\iota^{\alpha}} & TX \dar{T \alpha}  \\
      A \rar{0} & TA
    \end{tikzcd}
  \end{equation}
  and we write $p^{\alpha}=\iota^{\alpha}p$.
\end{definition}

\begin{lemma}\label{tangent-involution}
  The arrow
  \begin{equation*}
    X\tensor[_{\beta}]{\times}{_{\beta}}X\tensor[_{\alpha}]{\times}{_{\pi^{\alpha}}}T^{\alpha}X \xrightarrow{\xi} X\tensor[_{\beta}]{\times}{_{\beta}}X\tensor[_{\alpha}]{\times}{_{\pi^{\alpha}}}T^{\alpha}X \\
  \end{equation*}
  defined by $(\pi_1,\pi_0,(\pi_00,\pi_10,\pi_2\iota^{\alpha})T\mu)$ is an involution (i.e. $\xi\xi = 1_C$).
  \begin{proof}
    \begin{align}
      \xi\xi &= (\pi_1,\pi_0,(T\mu)(\pi_00,\pi_10,\pi_2))(\pi_1,\pi_0,(T\mu)(\pi_00,\pi_10,\pi_2))\\
      &=(\pi_0,\pi_1,(T\mu)(\pi_10,\pi_00,(T\mu)(\pi_00,\pi_10,\pi_2)))\\
      &=(\pi_0,\pi_1,(T\mu)(\pi_10,\pi_10,\pi_2))\\
      &=(\pi_0,\pi_1,\pi_2)
    \end{align}
  \end{proof}
\end{lemma}

\subsection{Splitting of the Vertical Bundle}\label{sec:splitting-of-the-vertical-bundle}

In this section we prove that when we assume the existence of an arrow $e:B \rightarrow X$ such that $\beta e \beta = \beta$ that we can split the bundle $T^{\alpha}_BX$ as $X \tensor[_{\beta}]{\times}{_{\pi^{\alpha}_M}}T^{\alpha}_M X$. The isomorphism exhibiting the splitting is the pullback of the involution $\xi$ defined in \Cref{tangent-involution} along the arrow $(\pi_0,\pi_1p^{\alpha}_B,\pi_1\iota_B):X \tensor[_{\beta}]{\times}{_{\pi^s_M}}T^s_M X \rightarrow X\tensor[_{\beta}]{\times}{_{\beta}}X\tensor[_{\alpha}]{\times}{_{\pi^{\alpha}}}T^{\alpha}X$.

\begin{notation}
  Let $e:B \rightarrow X$ be an arrow in $\mathbb X$ satisfying $\beta e \beta = \beta$. Note that if $\beta$ is an epimorphism then $e$ is a section of $\beta$.
\end{notation}

\begin{definition}\label{pullbacks-of-the-tangent-bundle}
  The arrow $\pi^{\alpha}_B:T^{\alpha}_B X \rightarrow B$ is defined by the pullback
  \begin{equation}
    \begin{tikzcd}
      T^{\alpha}_B X \rar{\iota_B} \dar{\pi^{\alpha}_B} & T^{\alpha} X \dar{\iota^{\alpha} p}\\
      B \rar{e} & X
    \end{tikzcd}
  \end{equation}
  and we write $p^{\alpha}_B = \iota_B p^{\alpha}$.
\end{definition}

\begin{lemma}\label{middle-square}
  The square
  \begin{equation}
    \begin{tikzcd}[column sep=2cm]
      X\tensor[_{\beta}]{\times}{_{\pi^{\alpha}_B}}T^{\alpha}_B X \rar{(\pi_0,\pi_1p^{\alpha}_B,\pi_1\iota_B)} \dar{\pi_0} & X\tensor[_{\beta}]{\times}{_{\beta}}X\tensor[_{\alpha}]{\times}{_{\pi^{\alpha}}}T^{\alpha}X \dar{(\pi_0,\pi_1,\pi_2 p^{\alpha})}\\
      X \rar{(1_X,\beta e,\beta e)} & X\tensor[_{\beta}]{\times}{_{\beta}}X\tensor[_{\alpha}]{\times}{_{\alpha}}X
    \end{tikzcd}
  \end{equation}
  is a pullback in $\mathbb{X}$.
  \begin{proof}
    First the right square of
    \begin{equation*}
      \begin{tikzcd}
        X \tensor[_{\beta}]{\times}{_{\pi^{\alpha}_B}} T^{\alpha}_B X \rar{\pi_1} \dar{\pi_0} & T^{\alpha}_B X \rar{\iota_B} \dar{\pi^{\alpha}_B} & T^{\alpha}X \dar{p^{\alpha}}\\
          X \rar{\beta} & B \rar{e} & X
      \end{tikzcd}
    \end{equation*}
    is a pullback by definition of $\pi^{\alpha}_B$ and the left square is immediately seen to be a pullback. Therefore the whole square is a pullback. Therefore the outer square of
    \begin{equation}
      \begin{tikzcd}[column sep=2cm]
        X\tensor[_{\beta}]{\times}{_{\pi^{\alpha}_B}}T^{\alpha}_B X \rar{(\pi_0,\pi_1p^{\alpha}_B p,\pi_1\iota_B)} \dar{\pi_0} & X\tensor[_{\beta}]{\times}{_{\beta}}X\tensor[_{\alpha}]{\times}{_{\pi^{\alpha}}}T^{\alpha}X \dar{(\pi_0,\pi_1,\pi_2 \iota^{\alpha}p)} \rar{\pi_2} & T^{\alpha}X \dar{p^{\alpha}}\\
        X \rar{(1_X,\beta e,\beta e)} & X\tensor[_{\beta}]{\times}{_{\beta}}X\tensor[_{\alpha}]{\times}{_{\alpha}}X \rar{\pi_2} & X
      \end{tikzcd}
    \end{equation}
    is a pullback by the definition of $T^{\alpha}_B X$ and the right hand square is immediately seen to be a pullback also. Therefore by the pullback pasting lemma the left hand square is a pullback as required.
  \end{proof}
\end{lemma}

\begin{lemma}\label{lower-square}
  The square
  \begin{equation}\label{involution-pullback}
    \begin{tikzcd}
      X \dar{(\beta e,1_X,1_X)} \rar{1_X} & X \dar{(1_X,\beta e,\beta e)} \\
      X\tensor[_{\beta}]{\times}{_{\beta}}X\tensor[_{\alpha}]{\times}{_{\alpha}}X \rar{\check{\xi}} & X\tensor[_{\beta}]{\times}{_{\beta}}X\tensor[_{\alpha}]{\times}{_{\alpha}}X
    \end{tikzcd}
  \end{equation}
  is a pullback.
  \begin{proof}
    First the square commutes because
    \begin{align}
      (\beta e , 1_X, 1_X)\check{\xi} &= (\beta e , 1_X, 1_X)(\pi_1,\pi_0,\pi_0\circ_{\pi_1}\pi_2)\\
      &=(1_X,\beta e,\beta e\circ_{1_X}1_X) = (1_X,\beta e,\beta e)
    \end{align}
    Second the diagram is a pullback because $1_x$ and $\check{\xi}$ are isomorphisms (and the diagram commutes).
  \end{proof}
\end{lemma}

\begin{lemma}\label{outer-square}
  The square
  \begin{equation}
    \begin{tikzcd}[column sep = 2cm]
      T^{\alpha}X \rar{(p^{\alpha}\beta e, p^{\alpha},1_{T^{\alpha}X})} \dar{p^{\alpha}} & X\tensor[_{\beta}]{\times}{_{\beta}}X\tensor[_{\alpha}]{\times}{_{\pi^{\alpha}}}T^{\alpha}X \dar{(\pi_0,\pi_1,\pi_2p^{\alpha})}\\
      X \rar{(\beta e, 1_X, 1_X)} & X\tensor[_{\beta}]{\times}{_{\beta}}X\tensor[_{\alpha}]{\times}{_{\alpha}}X
    \end{tikzcd}
  \end{equation}
  is a pullback.
  \begin{proof}
    In the diagram
    \begin{equation}
      \begin{tikzcd}[column sep = 2cm]
        T^{\alpha}X \rar{(p^{\alpha}\beta e, p^{\alpha},1_{T^{\alpha}X})} \dar{\iota^{\alpha}p} & X\tensor[_{\beta}]{\times}{_{\beta}}X\tensor[_{\alpha}]{\times}{_{\pi^{\alpha}}}T^{\alpha}X \dar{(\pi_0,\pi_1,\pi_2p^{\alpha})} \rar{\pi_2} & T^{\alpha}X \dar{p^{\alpha}}\\
        X \rar{(\beta e, 1_X, 1_X)} & X\tensor[_{\beta}]{\times}{_{\beta}}X\tensor[_{\alpha}]{\times}{_{\alpha}}X \rar{\pi_2} & X
      \end{tikzcd}
    \end{equation}
    the right hand square is immediately seen to be a pullback and the whole square is the trivial pullback. Therefore by the pasting lemma the left square is a pullback as required.
  \end{proof}
\end{lemma}

\begin{proposition}
  The objects $X\tensor[_{\beta}]{\times}{_{\pi^{\alpha}_B}}T^{\alpha}_B X$ and $T^{\alpha}X$ are isomorphic.
  \begin{proof}
    In the cube
    \begin{equation}
\begin{tikzcd}
T^{\alpha}X \drar \arrow{rrr}{(\iota^{\alpha}p\beta e, \iota^{\alpha}p,1_{T^{\alpha}X})} \arrow{ddd}{p^{\alpha}} & {} & {} & X\tensor[_{\beta}]{\times}{_{\beta}}X\tensor[_{\alpha}]{\times}{_{\pi^{\alpha}}}T^{\alpha}X \arrow{ddd}[swap]{(\pi_0,\pi_1,\pi_2p^{\alpha})} \dlar{\xi}\\
        {} & X\tensor[_{\beta}]{\times}{_{\pi^{\alpha}_B}}T^{\alpha}_B X \rar[yshift=1ex]{(\pi_0,\pi_1p^{\alpha}_B,\pi_1\iota_B)} \dar{\pi_0} & X\tensor[_{\beta}]{\times}{_{\beta}}X\tensor[_{\alpha}]{\times}{_{\pi^{\alpha}}}T^{\alpha}X \dar{(\pi_0,\pi_1,\pi_2 p^{\alpha})} & {}\\
        {} & X \rar{(1_X,\beta e,\beta e)} & X\tensor[_{\beta}]{\times}{_{\beta}}X\tensor[_{\alpha}]{\times}{_{\alpha}}X & {}\\
        X \urar{1_X} \arrow{rrr}{(\beta e, 1_X, 1_X)} & {} & {} & X\tensor[_{\beta}]{\times}{_{\beta}}X\tensor[_{\alpha}]{\times}{_{\alpha}}X \ular{\check{\xi}}
      \end{tikzcd}
    \end{equation}
    the middle square, the lower square and the outer square are pullbacks by \Cref{middle-square}, \Cref{lower-square} and \Cref{outer-square} respectively.
  \end{proof}
\end{proposition}

\bibliography{references}
\bibliographystyle{plain}

\end{document}